\documentclass[12pt]{article}
\usepackage{graphicx}
\usepackage{amsmath,amsthm,amssymb,enumerate}
\usepackage{euscript,mathrsfs}
\usepackage{color}
\usepackage{dsfont}
\usepackage{url}
\usepackage[left=2cm,right=2cm,top=3.5cm,bottom=3.5cm]{geometry}
\usepackage{color}
\usepackage[framemethod=tikz]{mdframed}
\allowdisplaybreaks
\usepackage{esint}

\usepackage{soul}

\catcode`\@=11 \@addtoreset{equation}{section}

\catcode`\@=12

\newtheorem{Theorem}{Theorem}[section]
\newtheorem{Proposition}[Theorem]{Proposition}
\newtheorem{Lemma}[Theorem]{Lemma}
\newtheorem{Corollary}[Theorem]{Corollary}

\theoremstyle{definition}
\newtheorem{Definition}[Theorem]{Definition}

\newtheorem{Remark}[Theorem]{Remark}

\newcommand{\bTheorem}[1]{
	\begin{Theorem} \label{T#1} }
	\newcommand{\eT}{\end{Theorem}}

\newcommand{\bProposition}[1]{
	\begin{Proposition} \label{P#1}}
	\newcommand{\eP}{\end{Proposition}}

\newcommand{\bLemma}[1]{
	\begin{Lemma} \label{L#1} }
	\newcommand{\eL}{\end{Lemma}}

\newcommand{\bCorollary}[1]{
	\begin{Corollary} \label{C#1} }
	\newcommand{\eC}{\end{Corollary}}

\newcommand{\bRemark}[1]{
	\begin{Remark} \label{R#1} }
	\newcommand{\eR}{\end{Remark}}

\newcommand{\bDefinition}[1]{
	\begin{Definition} \label{D#1} }
	\newcommand{\eD}{\end{Definition}}

\newcommand{\tvS}{\widetilde{S}}

\newcommand{\avintO}[1]{\fint_{\Omega} #1 \dx}

\newcommand{\tvm}{\widetilde{\vc{m}}}

\newcommand{\bfphi}{\boldsymbol{\varphi}}

\newcommand{\bFormula}[1]{
	\begin{equation} \label{#1}}
	\newcommand{\eF}{\end{equation}}

\newcommand{\Ov}[1]{\overline{#1}}

\newcommand{\vr}{\varrho}

\newcommand{\tvr}{\wtilde \vr}

\newcommand{\tvt}{\wtilde \vt}

\newcommand{\vt}{\vartheta}

\newcommand{\vm}{\vc{m}}

\newcommand{\vc}[1]{{\bf #1}}

\newcommand{\Div}{{\rm div}_x}
\newcommand{\Grad}{\nabla_x}

\newcommand{\dx}{\,{\rm d} {x}}

\newcommand{\dt}{\,{\rm d} t }

\newcommand{\intO}[1]{\int_{\Omega} #1 \ \dx}

\newcommand{\D}{{\rm d}}

\newcommand{\ep}{\varepsilon}

\newcommand{\br}{ \nonumber \\ }

\def\softd{{\leavevmode\setbox1=\hbox{d}%
		\hbox to 1.05\wd1{d\kern-0.4ex{\char039}\hss}}}
\definecolor{Cgrey}{rgb}{0.85,0.85,0.85}
\definecolor{Cblue}{rgb}{0.50,0.85,0.85}
\definecolor{Cred}{rgb}{1,0,0}
\definecolor{fancy}{rgb}{0.10,0.85,0.10}
\definecolor{amaranth}{rgb}{0.9, 0.17, 0.31}

\newcommand\Cbox[2]{%
	\newbox\contentbox%
	\newbox\bkgdbox%
	\setbox\contentbox\hbox to \hsize{%
		\vtop{
			\kern\columnsep
			\hbox to \hsize{%
				\kern\columnsep%
				\advance\hsize by -2\columnsep%
				\setlength{\textwidth}{\hsize}%
				\vbox{
					\parskip=\baselineskip
					\parindent=0bp
					#2
				}%
				\kern\columnsep%
			}%
			\kern\columnsep%
		}%
	}%
	\setbox\bkgdbox\vbox{
		\color{#1}
		\hrule width  \wd\contentbox %
		height \ht\contentbox %
		depth  \dp\contentbox
		\color{black}
	}%
	\wd\bkgdbox=0bp%
	\vbox{\hbox to \hsize{\box\bkgdbox\box\contentbox}}%
	\vskip\baselineskip%
}

\mdfdefinestyle{MyFrame}{%
	linecolor=black,
	outerlinewidth=1pt,
	roundcorner=5pt,
	innertopmargin=\baselineskip,
	innerbottommargin=\baselineskip,
	innerrightmargin=10pt,
	innerleftmargin=10pt,
	backgroundcolor=white!20!white}

\newcommand{\wtilde}{\widetilde}

\allowdisplaybreaks

\begin{document}


\title{\bf Statistical solutions to the Euler system of gas dynamics}

\author{Eduard Feireisl 	\thanks{
		The work of E.F.\ was partially supported by the
		Czech Sciences Foundation (GA\v CR), Grant Agreement
		24--11034S. The Institute of Mathematics of the Academy of Sciences of
		the Czech Republic is supported by RVO:67985840.
		E.F.\ is a member of the Ne\v cas Center for Mathematical Modelling.} }

\date{}

\maketitle

\medskip

\centerline{Institute of Mathematics of the Academy of Sciences of the Czech Republic}
\centerline{\v Zitn\' a 25, CZ-115 67 Praha 1, Czech Republic}
\centerline{feireisl@math.cas.cz}

\maketitle

\begin{abstract}
	
We consider the (complete) Euler system describing the motion of a compressible perfect fluid. We propose a platform suitable 
for constructing the statistical solutions. The main ingredients of our approach include: 
\begin{itemize}
\item The concept of dissipative (measure--valued) solution to the Euler system.

\item A single step selection procedure based on minimizing the Bregman divergence of a given solution to the maximal entropy equilibrium.

\item A construction of a Markov semigroup via push forward measures.

\end{itemize}	

\end{abstract}


{\small

\noindent

\medbreak


}

\bigskip

\section{Introduction}

The Euler system of gas dynamics belongs to the iconic examples of non--linear 
hyperbolic systems of conservation laws. The motion of a compressible perfect fluid (gas) confined to a (bounded) container 
$\Omega \subset R^d$, $d=1,2,3$ is described by means of three field equations 
expressing the basic underlying physical principles of 
\bigskip

\noindent 
{\bf the mass conservation}
\begin{equation} \label{i1}
	\partial_t \vr + \Div \vm = 0;
\end{equation}	
{\bf the momentum balance}
\begin{equation} \label{i2}
\partial_t \vm + \Div \left( \frac{\vm \otimes \vm}{\vr}  + p \mathbb{I} \right) = 0;
\end{equation}	
{\bf the energy balance}
\begin{equation} \label{i3}
\partial_t E + \Div \left( (E + p) \frac{\vm}{\vr} \right) = 0.
\end{equation}	
We retain the mass density $\vr = \vr(t,x)$ and the momentum 
$\vm = \vm(t,x)$ as the state variables, while the pressure 
$p$ as well as the energy $E$ will be given constitutively. Specifically, we write the energy in the form 
\begin{equation} \label{i4}
E = \frac{1}{2} \frac{|\vm|^2 }{\vr} + \vr e, 
\end{equation} 
where the internal energy $e$ is related to the pressure by the polytropic equation of state 
\begin{equation} \label{i5}
p = (\gamma - 1) \vr e, \ \gamma > 1.	
\end{equation}	
Finally, we introduce the entropy $s$ interrelated to the internal energy and the pressure through Gibbs' equation 
\begin{equation} \label{i6}
\vt D s = De + p D \left( \frac{1}{\vr} \right), 
\end{equation}
where the new variable $\vt$ is the absolute temperature. For the sake of simplicity, and actually without loss of generality, we may suppose 
\begin{equation} \label{i7} 
p = \vr \vt,\ e = c_v \vt,\ c_v = \frac{1}{\gamma - 1},\ 
s = c_v \log(\vt) - \log(\vr), 
\end{equation}
where $c_v$ is the specific heat at constant volume. Observe that the polytropic equation of state \eqref{i5} is incomplete (cf. the monograph 
of Benzoni--Gavage and Serre \cite{BenSer} ), and its specific form \eqref{i7} is determined by a proper choice of the absolute temperature $\vt$. 

It is convenient to complete the family phase variables $\vr$, $\vm$, by adding the total entropy $S(t,x) = \vr(t,x) s(\vr(t,x), \vt(t,x))$.
The problem is formally closed by imposing

\noindent
{\bf the impermeability boundary conditions}
\begin{equation} \label{i8}
\vm \cdot \vc{n}|_{\partial \Omega}; 
\end{equation}
and

\noindent
{\bf the initial conditions}
\begin{equation} \label{i9}
\vr(0, \cdot) = \vr_0,\ 
\vm(0, \cdot) = \vm_0,\ S = \vr s (0, \cdot) = S_0.	 
\end{equation}

As is well know, the Euler system is locally well posed in a class of sufficiently regular solutions, however, singularities, typically shock waves, may appear in a finite time no matter how smooth and small the initial data are cf. e.g. the monograph by Dafermos \cite{D4a}. The weak (distributional) solutions are able to accommodate singularities but are, in general, not unique cf. Smoller \cite{SMO}. 

To 
recover uniqueness, the entropy admissibility criterion has been proposed, namely

\noindent
{\bf the entropy inequality} 
\begin{equation} \label{i10}
\partial_t (\vr s) + \Div (s \vc{m}) \geq 0.
\end{equation}	
The inequality \eqref{i10} can be seen as a mathematical formulation of the Second law of thermodynamics, see Dafermos \cite{Daf4}.

There has been a common belief that the admissibility condition \eqref{i10} may recover uniqueness in the class of weak solutions. This is partially 
true in the one space dimension $d=1$. Besides the classical results that can be found in \cite{D4a}, there is a number of 
quite recent results by Vasseur et al. \cite{ChenVas}, \cite{KaVawa}    .
However, the recent adaptation of the method of convex integration initiated in the seminal work of De Lellis and Sz\' ekelyhidi 
\cite{DelSze2}, \cite{DelSze3} revealed a number rather disturbing facts concerning the well posedness of the Euler system in the class of admissible weak solutions in the space dimension $d=2,3$:
\begin{itemize}
	\item There is a large class of initial data for which the 
Euler system admits infinitely many admissible weak solutions, 
see \cite{FeKlKrMa}. In particular, a multidimensional extension 
of the classical Riemann problem is ill posed, see \cite{ABKlKrMaMa}, \cite{KlKrMaMa}, at least for certain initial data.

\item The set of the ``wild'' initial data that give rise to infinitely many entropy admissible weak solutions is dense in the $L^p$-topology of the phase space \cite{ChFe2023}.

\item The ill--posedness of the problem persists under stochastic perturbations, see \cite{ChFeFl2019}.
	
\end{itemize}	

We consider the Euler system in the framework of \emph{statistical solutions}. These were introduced by 
Foias \cite{Foias}, and Vishik and 
Fursikov \cite{VisFur} in the context of the incompressible Navier--Stokes system. The original theory has been 
later developed by \emph{Roger Temam} and his collaborators 
in a series of works \cite{FMRT},
\cite{FoRoTe3}, \cite{FoRoTe1}. Very roughly indeed, a statistical solution is a stochastic process supported by the solution set of a given problem. The initial data are understood as random variables ranging in the associated phase space. 

In \cite{FanFei}, a direct construction of a statistical solution was proposed being simply a pushforward measure of the distribution of the initial data. In such a way, a Markov semigroup on the space of probability measures is obtained characterizing the distribution of the solution at any positive time $t > 0$. Obviously, such an approach requires the problem 
to admit a Borel measurable semigroup of solutions. 

The method of Krylov \cite{KrylNV} and its adaptation by Cardona and Kapitanski \cite{CorKap} was used in \cite{FanFei} 
to select a solution semigroup for the compressible Navier-Stokes system (cf. also Basari\' c \cite{Basa1}). Motivated 
by \cite{FeLmJu2025}, \cite{FeiLukYu}, we extend the class of 
admissible (weak) solutions of the Euler system to the dissipative 
measure valued solutions introduced in \cite{BreFei17}, \cite{BreFei17A}. In contrast with \cite{FanFei}, where an infinite family of cost functionals is needed to select a unique solution, we propose a single step selection criterion 
enforced by the Second law of thermodynamics. Specifically, 
we elaborate on the classical statement of Clausius (Rudolf Clausius, Poggendorff’s Annals of Physics 1865 (125), 400) considered as the most flagrant formulation of the Second law of thermodynamics: 

\begin{quotation}
	
	\emph{The energy of the world (closed system)  is constant; its entropy tends to a maximum.}
	
\end{quotation}
	
\noindent
Accordingly, we propose a selection criterion based on minimizing the distance to the equilibrium state that maximizes the entropy. 	
Once a proper selection of the solution mapping is made, the construction of the associated Markov semigroup is a routine matter. 

The paper is organized as follows. In Section \ref{mv}
we recall the definition of dissipative (measure--valued) solutions of the Euler system and their basic properties. 
In Section \ref{s}, we propose a single step selection criterion based on minimization of the distance to the maximal entropy equilibrium. In particular, we show that the energy defect of the selected solution vanishes for $t \to \infty$, which is a result of independent interest. The main result on the existence of a suitable statistical solution is stated in Section \ref{m}. 
The proof of of the basic properties of statistical solutions is given in Section \ref{p}. The paper is concluded by a short discussion in Section \ref{c}.

\section{Dissipative solutions}
\label{mv} 

Despite the vast set of the (wild) initial data, for which the Euler system admits infinitely many admissible \emph{weak} solutions, the \emph{existence} of weak solutions for general initial data remains 
an outstanding open problem. Following \cite{BreFei17} we extend the class
of solutions incorporating all possible limits of their consistent 
approximations, cf. \cite[Part II, Chapter 5]{FeLMMiSh}. We point out that, in contrast with the ``pure'' measure--valued approach 
proposed by e.g. by Fjordholm et al. \cite{FjKaMiTa}, \cite{FjMiTa1}, Gwiazda et al. \cite{GSWW}, Klingenberg et al. \cite{KlMaWi}, the dissipative solutions are represented by 
\emph{functions} rather than parametrized (Young) measures.

First, we introduce the kinetic energy
\begin{equation} \label{s1}
 	\frac{1}{2} \frac{|\vm|^2}{\vr} = \left\{ \begin{array}{l}
 		\frac{1}{2} \frac{|\vm |^2}{\vr} \ \mbox{if}\ \vr > 0,\\
 		0 \ \mbox{if}\ \vr = 0, \vm = 0, \\
 		\infty \ \mbox{otherwise.} \end{array}
 	\right.
\end{equation}  
Similarly, the internal energy and the pressure are defined as 
\begin{equation} \label{s2}
	p(\vr, S) = (\gamma - 1) \vr e (\vr, S) = \left\{ \begin{array}{l}
		\vr^\gamma \exp \left( \frac{S}{c_v \vr} \right) \ \mbox{if} \ \vr > 0, \\
		0 \ \mbox{if}\ \vr = 0,\ S \leq 0, \\
		\infty \ \mbox{otherwise}.
	\end{array}	
	\right.
\end{equation}	
Finally, the total energy reads 
\begin{equation} \label{mv3}
E(\vr, \vm, S) = \frac{1}{2} \frac{|\vm|^2}{\vr} + \vr e (\vr, S), 
\end{equation}	
where we have adopted \eqref{s1}, \eqref{s2}.
In accordance with \eqref{s1}, \eqref{s2}, the total energy $E(\vr, \vm, S)$ is a convex l.s.c. function of $(\vr, \vm, S) \in R^{d+2}$, strictly convex on its domain. 

The dissipative solutions introduced below will belong to the phase space 
\[
\mathcal{P} = \left\{ (\vr, \vm, S, \mathcal{E}_0) \ \Big| 
\intO{ E(\vr, \vm, S) } \leq \mathcal{E}_0,\ S \geq \underline{s} \vr,\ \intO{ \vr } \geq m_0 > 0 \right\},
\]
where $\underline{s} \in R$, $m_0 > 0$ are a given constants. It is easy to check that the set $\mathcal{P}$ is convex.

\begin{Definition} [\bf Dissipative solution] \label{mvD1}
	
We say that $(\vr, \vm, S, \mathcal{E}_0)$ is \emph{dissipative 
solution} of the Euler system \eqref{i1}--\eqref{i3},  
\eqref{i10}, with the boundary condition \eqref{i8}, and the initial data 
\[
(\vr_0, \vm_0, S_0, \mathcal{E}_0) \in \mathcal{P} 
\]
in the space time cylinder $(0, \infty) \times \Omega$
if there exists a parametrized measure 
\begin{equation} \label{mv4}
\mathcal{V}_{t,x} ,\ t \in (0, \infty),\ x \in \Omega,\ 
(t,x) \mapsto \mathcal{V}_{t,x} \in L^\infty_{\rm weak-*} ((0, \infty) \times \Omega; 
\mathfrak{P}(R^{d + 2})) \ \mbox{for a.a}\ t \in (0, \infty),\ 
x \in \Omega, 
\end{equation} 
and a matrix--valued measure 
\begin{equation} \label{mv5}
\mathfrak{C} \in L^\infty_{\rm weak-*}(0, \infty; \mathcal{M}^+ 
(\Ov{\Omega}; R^{d + 2})) 
\end{equation} 
such that the following holds:	
\begin{enumerate}
\item The trio $(\vr, \vm, S)(t,x)$ coincides with the expected 
value of the parametrized measure $\mathcal{V}_{t,x}$, 
\[
(\vr, \vm, S)(t,x) = \int_{R^{d+2}}(\tvr, \tvm, \tvS) \ \D 
\mathcal{V}_{t,x} \ \mbox{for a.a.} \ t \in (0, \infty),\ 
x \in \Omega;
\]	
\begin{align}
\intO{ E \Big( (\vr, \vm, S)(t \pm , \cdot) \Big) } 
&\leq \mathcal{E}_0,\ S(t \pm, \cdot) \geq \underline{s} 
\vr(t, \cdot),\ \intO{ \vr(t, \cdot) } \geq m_0  \ \mbox{for any} \ t > 0, \br 
\intO{ E \Big( (\vr, \vm, S)(0+ , \cdot) } 
&\leq \mathcal{E}_0,\ S(0+, \cdot) \geq \underline{s} 
\vr_0.
\nonumber
\end{align}

\item The equation of continuity \eqref{i1} is satisfied in the sense of distributions, 
\begin{equation} \label{mv6}
\int_0^\infty \intO{ \Big[ \vr \partial_t \varphi + 
	\vm \cdot \Grad \varphi \Big] } \dt = \intO{ \vr_0 \varphi (0, \cdot)} 
\end{equation}
for any $\varphi \in C^1_c([0, \infty) \times \Ov{\Omega})$. 

\item The momentum equation \ref{i2} is replaced by the integral identity 
\begin{align} 
\int_0^\infty &\intO{ \left[ 
\vm \cdot \partial_t \bfphi + \int_{R^{d+2}}   \mathds{1}_{\tvr > 0} \left( \frac{\tvm \otimes \tvm}{\tvr} + 
p(\tvr, \tvS) \right) \D \mathcal{V}_{t,x} \cdot \Grad \bfphi   \right] } \dt \br 
&= - \int_0^\infty \left( \int_{\Ov{\Omega}} \D \mathfrak{C} : \Grad \bfphi \right) \dt - \intO{ \vm_0 \cdot \bfphi (0, \cdot) }
\label{mv7}
\end{align}	
for any $\bfphi \in C^1_c([0, \infty) \times \Ov{\Omega}; R^d)$, $\bfphi \cdot \vc{n}|_{\partial \Omega} = 0$.	

\item 
The entropy balance \eqref{i10} holds as an integral inequality 
\begin{equation} \label{mv8}
\int_0^\infty \intO{ \left[ S \partial_t \varphi + \int_{R^{d+2}} 
\left( \mathds{1}_{\tvr > 0} \tvS \frac{\tvm}{\tvr} \right) \D \mathcal{V}_{t,x} \cdot \Grad \varphi \right] } \dt \leq - 
\intO{ S_0 \varphi (0, \cdot) }		
\end{equation}	
for any $\varphi \in C^1_c([0, \infty) \times \Ov{\Omega})$, $\varphi \geq 0$. 

\item The energy balance \eqref{i3} is replaced by the total energy balance 
\begin{equation} \label{mv9}
\intO{ \int_{R^{d + 2}} E( \tvr, \tvm, \tvS) \D \mathcal{V}_{\tau, x} } + c(d, \gamma) \int_{\Ov{\Omega}} 
\D {\rm trace} [\mathfrak{C} (\tau, \cdot) ]  \leq \mathcal{E}_0,\ c(d, \gamma) = \min  
\end{equation} 	
for a.a. $\tau \geq 0$.

\end{enumerate}	

\end{Definition}	

We refer the reader to \cite[Chapter 5]{FeLMMiSh} for 
a detailed discussion of the concept of dissipative solutions to various problems in fluid mechanics and their use in numerical analysis.
They have been identified as limits of the so-called consistent approximation. The parametric family $(\mathcal{V}_{t,x})_{t \in (0, \infty),x \in \Omega}$ coincides 
with the Young measure characterizing the oscillations of the approximate sequence, while $\mathfrak{C}$ represents the oscillation defect measure.

We record the following properties of dissipative solutions:

\begin{itemize}
\item {\bf Global existence.} The Euler system admits global--in--time dissipative solutions for any choice of the initial data $(\vr_0, \vm_0, S_0, \mathcal{E}_0) \in \mathcal{P}$, see \cite[Proposition 3.8]{BreFeiHof19C}.

\item {\bf Weak continuity.} 
The functions $\vr$, $\vm$ belong to the class 
\begin{equation} \label{mv10}
\vr,\ \vm \in BC_{\rm weak}([0, \infty); L^\gamma(\Omega)); 
\end{equation}	
the total entropy admits one-sided limits 
\begin{equation} \label{mv11}
\lim_{t \to \tau \pm} \intO{ S (t, \cdot) \phi  } ,\
\lim_{t \to 0+} \intO{ S (t, \cdot) \phi} 
\ \mbox{for any} \ \phi \in L^{\gamma'}(\Omega). 
\end{equation}	
In particular $S(t\pm)$ can be identified with a function in $L^\gamma(\Omega)$, where we set $S(0-) = S_0$. There holds 
\begin{equation} \label{mv12}
S(t-) \leq S(t+) \ \mbox{for any}\ t \geq 0,\ 
\| S(t\pm) \|_{L^\gamma(\Omega)} \leq \Ov{S} \ \mbox{for all}\ 
t \geq 0,
\end{equation} 
see \cite[Chapter 5, Section 5.2]{FeLMMiSh}. 

\item {\bf Compatibility.} If 
$(\vr, \vm, S) \in C^1([ T_1, T_2] \times \Ov{\Omega}; R^d)$, 
$\inf_{[T_1, T_2] \times \Ov{\Omega}} \vr > 0$, and 
\[
\intO{ E(\vr, \vm, S) (t, \cdot) } = \mathcal{E}_0 \ 
\mbox{for a.a.}\ t \in (T_1, T_2), 
\]
then $(\vr, \vm, S)$ is a classical solution of the Euler system 
in $[T_1, T_2]$, see \cite[Chapter 5, Theorem 5.7]{FeLMMiSh}.

\item{\bf Weak--strong uniqueness.} 
The dissipative solution coincides with the strong $W^{1,\infty}$ solution of the Euler system as long as the latter solution exists, \cite[Chapter 6, Theorem 6.2]{FeLMMiSh}.
 
\end{itemize}

\section{Selection criterion}
\label{s}

Similarly to their weak counterparts, the dissipative solutions are not unique for a vast class of physically admissible initial data. In \cite{BreFeiHof19C}, a selection procedure in the spirit of Krylov \cite{KrylNV} and Cardona, Kapitanski \cite{CorKap} was proposed to identify a unique physically relevant solution. This process consists of a successive minimization of a countable family of ``cost functionals'' separating points in the phase space $\mathcal{P}$. A proper choice of these as well as the order of selection is, however, arbitrary and any ``physical'' justification is unclear. A significant progress has been made in \cite{FeiLM2025II}, where the selection process was reduced to merely two steps, and even a single step selection was proposed. Developing the ideas of \cite{FeiLM2025II}, we consider a single step selection process minimizing the Bregman distance 
of solutions to the constant equilibrium maximizing the entropy. 

\subsection{Stable equilibrium state} 
It follows directly from the equation of continuity \eqref{mv6} that the total mass 
\begin{equation} \label{s3}
	\intO{ \vr(t, \cdot) } =  \intO{\vr_0 } = M_0 > 0 
\end{equation}
is conserved. Given the total energy $\mathcal{E}_0$, the stable 
equilibrium density-temperature is given by the constant state 
$(\Ov{\vr}, \Ov{\vt})$,   
\[
\Ov{\vr} = \avintO{ \vr_0 },\ \mathcal{E}_0 = c_v \intO{ 
\Ov{\vr} \Ov{\vt} } \ \Rightarrow \ \Ov{\vt} = \frac{\mathcal{E}_0} {c_v M_0}.
\]

Next, we introduce the Bregman distance (divergence) associated to the convex energy $E$, 
\begin{align}
E &\left( \vr, \vm, S \Big| \tvr, \tvm, \tvS \right) \br
&= E( \vr, \vm, S) - \frac{\partial E (\tvr, \tvm, \tvS) }{\partial \vr} (\vr - \tvr) -  \frac{\partial E (\tvr, \tvm, \tvS) }{\partial \vm} \cdot (\vm - \tvm) -  \frac{\partial E (\tvr, \tvm, \tvS) }{\partial S} (S - \tvS) \br &\quad - E(\tvr, \tvm, \tvS).
\nonumber
\end{align}
We refer to \cite[Chapter 5]{FeLMMiSh} for a detailed discussion of the concept of Bregman divergence associated to the total energy $E$. 

Evaluating the Bregman distance of a dissipative solution $(\vr, \vm, S)$ to the equilibrium solution 
$(\Ov{\vr}, 0 , \Ov{\vr} s(\Ov{\vr}, \Ov{\vt}))$, and integrating over $\Omega$, we obtain 
\begin{equation} \label{ss1}
\intO{ E \left( \vr, \vm, S \Big| \Ov{\vr} , 0 , \Ov{\vr} s(\Ov{\vr}, \Ov{\vt}) \right) } = \intO{ \Big[  E(\vr, \vm, S) - \Ov{\vt} S ] } + \Ov{\vt} M_0 s(\Ov{\vr}, \Ov{\vt}) - \mathcal{E}_0,
\end{equation}	
cf. \cite[Chapter 5]{FeLMMiSh}. 
This motivates the choice of the cost functional 
\begin{equation} \label{s8}
	\mathcal{F}(\vr, \vm, S) = 
	\intO{ \Big[  E(\vr, \vm, S) - \Ov{\vt} S \Big] }, 
	\Ov{\vt} = \frac{\mathcal{E}_0}{c_v M_0}
\end{equation}
considered already in \cite{FeiLM2025II}.

\subsection{A single step selection}

Let 
\begin{align}
	\mathcal{U}[\vr_0, \vm_0, S_0, \mathcal{E}_0] = 
	&\left\{ (\vr, \vm, S) \ \Big|\ (\vr, \vm, S) 
	\ \mbox{is a dissipative solution in}\ (0,\infty) \times 
	\Omega \right. \br
	&\mbox{emanating from the data}\ 
	(\vr_0, \vm_0, S_0; \mathcal{E}_0) \Big\}
\nonumber	
\end{align}
be the set of all dissipative solutions associated to the initial data $(\vr_0, \vm_0, S_0; \mathcal{E}_0)$.
We introduce the weighted time-measure 
\[
\omega \dt = \exp(-t) \dt, 
\]
along with the associated weighted spaces 
\[
L^q_{\omega}(0, \infty; L^r(\Omega; R^{d+2})), \ L^q_{\omega}(0, \infty; W^{-\ell,2}(\Omega; R^{d+2})), \ \ell > d.
\] 
We refer to \cite[Chapter 5, Theorem 5.2]{FeLMMiSh} for the proof of the following properties of the solution set 
$\mathcal{U}[\vr_0, \vm_0, S_0, \mathcal{E}_0]$:
\begin{enumerate}
	\item The set $\mathcal{U}[\vr_0, \vm_0, S_0, \mathcal{E}_0]$ is non--empty and convex whenever $(\vr_0, \vm_0, S_0, \mathcal{E}_0) \in \mathcal{P}$.
	\item The set $\mathcal{U}[\vr_0, \vm_0, S_0, \mathcal{E}_0]$ is closed and bounded in the Banach space
	\[
	L^q_{\omega}(0, \infty; L^r(\Omega; R^{d+2})),\ 1 \leq q < \infty,\ 1 \leq r \leq \frac{2\gamma}{\gamma + 1}.
	\]
	\item The set $\mathcal{U}[\vr_0, \vm_0, S_0, \mathcal{E}_0]$ is compact in the Banach space 
	\[
	L^q_{\omega}(0, \infty; W^{-\ell,2}(\Omega; R^{d+2})),\ \ell > d.
	\]
\end{enumerate}	

We are ready to introduce the concept of \emph{admissible dissipative solution} to the Euler system selected in a single step. 

\begin{Definition}[{\bf Admissible dissipative solution}] \label{sD1}
We say that $(\vr, \vm, S, \mathcal{E}_0)$ is \emph{admissible dissipative solution} to the Euler system with the initial data 
\[
(\vr_0, \vm_0, S_0, \mathcal{E}_0) \in \mathcal{P} 
\]
if 
\begin{equation} \label{ss2}
(\vr, \vm, S) = {\rm arg min}_{\mathcal{U}[\vr_0, \vm_0, S_0, \mathcal{E}_0]} \int_0^\infty 
\exp(-t) \left( \mathcal{F} (\tvr, \tvm, \tvS) (t, \cdot) \right) \dt, 
\end{equation}	
meaning $(\vr, \vm, S) \in \mathcal{U}[\vr_0, \vm_0, S_0, \mathcal{E}_0]$ and 
\begin{align} 
\int_0^\infty & \exp(-t) \left( \intO{ \Big[  E(\vr, \vm, S) - \Ov{\vt} S \Big] } \right) \dt \leq 
\int_0^\infty & \exp(-t) \left( \intO{ \Big[  E(\tvr, \tvm, \tvS) - \Ov{\vt} \tvS \Big] } \right) \dt
\label{ss3}	
\end{align}	
for all $(\tvr, \tvm, \tvS) \in \mathcal{U}[\vr_0, \vm_0, S_0, \mathcal{E}_0]$.
\end{Definition}

As shown in \cite{FeiLM2025II}, admissible dissipative solution enjoy the following properties:
\begin{enumerate}
\item {\bf Existence and uniqueness.} 
For any initial data $(\vr_0, \vm_0, S_0, \mathcal{E}_0)$, the Euler system admits a unique admissible dissipative solution.  
\item {\bf Semigroup property.} The solution mapping 
\[
\mathcal{S}: (\vr_0, \vm_0, S_0, \mathcal{E}_0) \mapsto (\vr(t, \cdot), \vm(t, \cdot), S(t-, \cdot), \mathcal{E}_0))
\]
enjoys the semigroup property:
\begin{align}
\mathcal{S} [\vr_0, \vm_0, S_0, \mathcal{E}_0] (t+s, \cdot) = 	
\mathcal{S} \Big[ \mathcal{S}[\vr(t, \cdot, \vm(t, \cdot), S(t-,\cdot), \mathcal{E} \Big] (s, \cdot),\ t,s \geq 0.
\label{ss4}	
\end{align}	
 
\item {\bf Borel measurability.}	
The mapping 
\[
(\vr_0, \vm_0, S_0, \mathcal{E}_0) \in \mathcal{P} \mapsto 
(\vr(t, \cdot), \vm(t, \cdot), S(t-, \cdot), \mathcal{E}_0) \in \mathcal{P}
\]
is Borel measurable with respect to the $W^{-\ell,2}(\Omega; R^{d + 2}) \times R$ Hilbert topology for any $t \geq 0$.
\end{enumerate}

\subsection{Vanishing energy defect for $t \to \infty$}

In addition to the properties of admissible dissipative solutions proved in \cite{FeiLM2025II}, we show that 
any such solution is ``almost weak'' for $t \to \infty$. To this end, we first introduce the \emph{energy defect}
\begin{equation} \label{s7}
	d_E(t) = \mathcal{E}_0 - \intO{ E \Big( (\vr, \vm, S)(t+, \cdot) \Big) }. 
\end{equation}	
Note that, in accordance with \eqref{mv9}, 
\begin{align} 
d_E(t) &= \mathcal{E}_0 - \intO{ E \Big( (\vr, \vm, S)(t+, \cdot) \Big) } \br 
&\geq \intO{ \int_{R^{d + 2}} E( \tvr, \tvm, \tvS) \D \mathcal{V}_{t, x} } - 
\intO{ E \left( \int_{R^{d+2}} ( \tvr, \tvm, \tvS) \D \mathcal{V}_{t, x}               \right) } \br 
& + c(d, \gamma) \int_{\Ov{\Omega}} 
\D {\rm trace} [\mathfrak{C} (t, \cdot) ] 
\label{ss7} 
\end{align}	
for a.a. $t > 0$. By virtue of Jensen's inequality, the first term 
\[
\intO{ \int_{R^{d + 2}} E( \tvr, \tvm, \tvS) \D \mathcal{V}_{t, x} } - 
\intO{ E \left( \int_{R^{d+2}} ( \tvr, \tvm, \tvS) \D \mathcal{V}_{t, x} \right) }
\]
is always non-negative and vanishes only if $\mathcal{V}_{t, x}$ is the Dirac mass centred at $( \tvr, \tvm, \tvS)(t,x)$ for a.a. $x \in \Omega$.
Consequently, if $d_E = 0$, the admissible dissipative solution is a weak solution satisfying the Euler system in the following sense:  
\begin{align} 
\partial_t \vr + \Div \vm &= 0, \br
\partial_t \vm + \Div \left(  \frac{\tvm \otimes \tvm}{\tvr} + p(\vr, S) \mathbb{I} \right) &= 0, \br 
\partial_t S + \Div \left( S \frac{\vm}{\vr} \right) &\geq 0, \br
\intO{ E(\vr, \vm, S)(t, \cdot) } = \mathcal{E}_0 \ \mbox{for a.a.}\ t > 0,
\label{ss8}
\end{align}	
where the derivatives are understood in the sense of distributions. Note that the energy defect $d_E$ is related 
to the magnitude of the so-called Reynolds stress in the theory of deterministic turbulence, cf. Bardos et al. \cite{BaGhKa}.

Our objective is to show the following result that is of independent interest. 

\begin{Proposition}[{\bf Vanishing energy defect}] \label{sP1} Let $(\vr, \vm, S, \mathcal{E}_0)$ be an admissible 
dissipative solution of the Euler system in $(0,T) \times \Omega$. 

Then 	
\[
d_E(t) \equiv \mathcal{E}_0 - \intO{ E \Big( (\vr, \vm, S)(t+, \cdot) \Big) } \to 0 \ \mbox{as}\ t \to \infty.
\]
\end{Proposition}

The conclusion of Proposition \ref{sP1} can be interpreted in the way that any admissible dissipative solutions approaches a weak solution of the Euler system in the long run. The rest of this section is devoted to the proof of Proposition \ref{sP1}. We proceed in several steps.

\subsubsection{Preliminaries}
We recall that dissipative solutions
always satisfy
\begin{equation} \label{s4}
\intO{ E \Big( (\vr, \vm, S)(t\pm, \cdot) \Big) } \leq \mathcal{E}_0,
\ (\vr, \vm, S)(0-, \cdot) \equiv (\vr_0, \vm_0, S_0),
\end{equation}		
\begin{equation} \label{s5}
	S(t \pm, \cdot) \geq \underline{s} \vr.
\end{equation}	
Moreover, 
it follows from \eqref{s4} and \eqref{s5} that 
\begin{align} 
	\vr(t,x) &= 0 \ \Rightarrow \ \vm(t,x) = 0 \ \mbox{for a.a.}\ x \in \Omega, \br 
	\vr(t,x) &= 0 \ \Rightarrow \ S(t\pm,x) = 0 \ \mbox{for a.a.}\ x \in \Omega
\label{s6}
\end{align}

\subsubsection{Positive energy defect and jumps of the cost functional $\mathcal{F}$}
\label{ed}
Suppose that
\begin{equation} \label{s9} 
d_E(\tau) > 0 \ \mbox{for some}\ \tau > 0.
\end{equation}
Our aim is to show that the solution $(\vr, \vm, S, \mathcal{E}_0)$ can be extended beyond 
the time $\tau$ in such a way that the cost functional $\mathcal{F}$ evaluated at $\tau$ experiences a negative jump proportional 
to $d_E(\tau)$. 

First, using \eqref{s6}, we evaluate the temperature 
$\vt(\tau+, \cdot)$ from the relation 
\[
S(\tau+, \cdot) = \mathds{1}_{\vr(\tau, \cdot) > 0} \vr(\tau, \cdot) \Big( c_v \log(\vt(\tau+, \cdot)) - \log(\vr(\tau, \cdot)) \Big).
\]
Next, we introduce 
\[
\tvt (\tau, \cdot) = (1 + \ep) \vt(\tau+, \cdot) \ \mbox{whenever}\ \vr(\tau, \cdot) > 0  
\]
for $\ep > 0$ to be determined below. Denoting the internal energy 
\[
E_{\rm int}(\tau+, \cdot) = c_v \vr(\tau, \cdot) \vt(\tau+, \cdot), 
\]
we get
\begin{equation} \label{s10}
\widetilde{E}_{\rm int}(\tau, \cdot) \equiv c_v \vr(\tau, \cdot) \tvt(\tau, \cdot) = E_{\rm int}(\tau+, \cdot) + \ep c_v \vr(\tau, \cdot) \vt(\tau+, \cdot), 
\end{equation}
and
\begin{equation} \label{s11}
\widetilde{S}(\tau, \cdot) \equiv c_v \vr(\tau, \cdot) \log( \tvt(\tau, \cdot) ) - \vr(\tau, \cdot) \log( \vr (\tau, \cdot) ) 
= S(\tau+, \cdot)  + c_v \vr(\tau, \cdot) \log(1 + \ep)
\end{equation}

Now, we adjust $\ep > 0$ so that 
\begin{equation} \label{s12}
\widetilde{d}_E(\tau) = 
\mathcal{E}_0 - \intO{ E \Big( (\vr, \vm, \widetilde{S})(\tau, \cdot) \Big) } = 0, 
\end{equation}		
meaning 
\begin{equation} \label{s13}
\ep c_v \intO{ \vr(\tau, \cdot) \vt(\tau+, \cdot) } = d_E(\tau)
\ \Rightarrow \ 
\ep = \frac{d_E(\tau)}{\intO{ E_{\rm int}(\tau+, \cdot)}}. 
\end{equation} 

The relations \eqref{s12}, \eqref{s13} imply that we may continue the original solution $(\vr, \vm, S, \mathcal{E}_0)$ by concatenating it at the time $\tau$ with 
a solution $(\tvr, \tvm, \tvS, \mathcal{E})$ starting from the data  
\[
(\vr(\tau, \cdot), \vm(\tau, \cdot), \tvS(\tau, \cdot), \mathcal{E}_0) \ \mbox{for}\ t \geq \tau, 
\]
with the energy defect $\widetilde{d}_E(\tau) = 0$.

Next, we evaluate the cost functional at the time $\tau$ of the concatenated solution:
\begin{align} 
\mathcal{F} \Big( (\vr, \vm, \widetilde{S}) (\tau, \cdot) \Big) &=
\mathcal{F} \Big( (\vr, \vm, {S}) (\tau+, \cdot) \Big) \br  
&+ \ep \intO{ E_{\rm int}(\tau+, \cdot) } - 
\Ov{\vt} c_v \intO{ \vr(\tau, \cdot) \log(1 + \ep) }
\label{s14}
\end{align}
Moreover, using \eqref{s8}, \eqref{s13} we have
\begin{align}
\ep \intO{ E_{\rm int}(\tau+, \cdot) } &- 
\Ov{\vt} c_v \intO{ \vr(\tau, \cdot) \log(1 + \ep) } \br 
&= d_E(\tau) - \mathcal{E}_0 \log \left( 1 + \frac{d_E(\tau)}{\intO{ E_{\rm int}(\tau+, \cdot)}} \right).
\label{s15}
\end{align}	

Finally, we observe 
\[
\intO{ E_{\rm int}(\tau+, \cdot)} \leq 
\intO{ E\Big(\vr, \vm, S)(\tau+, \cdot) \Big) } = \mathcal{E}_0 - d_E(\tau);  
\]
whence
\begin{equation} \label{s17}
d_E(\tau) - \mathcal{E}_0 \log \left( 1 + \frac{d(\tau)}{\intO{ E_{\rm int}(\tau+, \cdot)}} \right)
\leq d_E(\tau) - \mathcal{E}_0 \log \left( 1 + \frac{d(\tau)}{\mathcal{E}_0 - d(\tau) } \right).
\end{equation}
Introducing the function
\begin{equation} \label{s18}
G: y \mapsto \mathcal{E}_0 \log \left( 1 + \frac{y}{\mathcal{E}_0 - y } \right) - y
= \mathcal{E}_0 \log \left( \frac{\mathcal{E}_0 }{\mathcal{E}_0 - y} \right) - y = 
\mathcal{E}_0 \log( \mathcal{E}_0 ) - \mathcal{E}_0 
\log(\mathcal{E}_0 - y) - y
\end{equation}
satisfies
\begin{equation} \label{s19}
G(0) = 0,\ G'(y) = \frac{\mathcal{E}_0}{\mathcal{E}_0 - y}  - 1	> 0 \ \mbox{as long as}\ y > 0.
\end{equation}
Thus we may infer
\begin{equation} \label{s20} 
\widetilde{d}_E   (\tau) = 0,\ 
\mathcal{F} \Big( (\vr, \vm, \widetilde{S}) (\tau, \cdot) \Big)
\leq \mathcal{F} \Big( (\vr, \vm, {S}) (\tau+, \cdot) \Big) -
G(d_e(\tau))
\end{equation}
for the concatenated solution, 
where $G$ is defined in \eqref{s18}.

\subsubsection{Vanishing defect for the selected solution as $t \to \infty$}

Having prepared the preliminary material, we are ready to prove Proposition \ref{sP1}.

Assume that $(\vr, \vm, S, \mathcal{E}_0)$ is an admissible dissipative solution, meaning the global minimizer of the cost functional 
\begin{equation} \label{s21}
\int_0^\infty \exp(-t) \left( \intO{ \mathcal{F} \Big( \vr, \vm, S)(t, \cdot) \Big)                } \right) \dt	
\end{equation}
on the set $\mathcal{U}(\vr_0, \vm_0, S_0, \mathcal{E}_0)$.
Our aim is to show that the associated defect vanishes in the long run, 
\begin{equation} \label{s22}
d_E(t) \to 0 \ \mbox{as}\ t \to \infty. 
\end{equation}

We show \eqref{s22} arguing by contradiction. If \eqref{s22} is not true, than there exists $\Ov{d} > 0$ such that 
\begin{align} 
	\limsup_{t \to \infty} d_E(t) &= \Ov{d} > 0, \label{s23}\\ 
	d_E(\tau_n) &\to \Ov{d} \ \mbox{for a certain sequence}\ 
	\tau_n \to \infty. 
	\label{s24}
\end{align}	

\begin{enumerate}
	\item
Since 
\[
d_E(\tau_n) = \mathcal{E}_0 - \intO{ E \Big( (\vr, \vm, S) (\tau_n+ ) \Big) },	
\]
we may rephrase \eqref{s23}, \eqref{s24} as 
\begin{align} 
	\liminf_{t \to \infty} \intO{E \Big( (\vr, \vm, S) (t+, \cdot ) \Big) }  &\geq \mathcal{E}_0 - \Ov{d} , \label{s25}\\ 
		\intO{E \Big( (\vr, \vm, S) (\tau_n +, \cdot ) \Big) } &\to \mathcal{E}_0 - \Ov{d} \ \mbox{for a certain sequence}\ 
	\tau_n \to \infty. 
	\label{s26}
\end{align}	
It follows from \eqref{s25}, \eqref{s26} that 
for any $\omega > 0$, there exists $n = n(\omega)$ such that 
\begin{equation} \label{s27}
\intO{E \Big( (\vr, \vm, S) (t+, \cdot ) \Big) } 
\geq \intO{E \Big( (\vr, \vm, S) (\tau_n +, \cdot ) \Big) } - \omega \ \mbox{for any}\ t \geq \tau_n,\ n \geq n(\omega).
\end{equation}	

\item As the total entropy admits a limit for $t \to \infty$, 
\[
\intO{ S (t\pm) } \nearrow \mathcal{S}_\infty < \infty 
\ \mbox{as}\ t \to \infty,
\]
we conclude, similarly to \eqref{s27},
\begin{align} 
&\intO{E \Big( (\vr, \vm, S) (t+, \cdot ) \Big) } - \Ov{\vt} 
\intO{ S(t+, \cdot)}\br &\quad \geq 	\intO{E \Big( (\vr, \vm, S) (\tau_n +, \cdot ) \Big) } - \Ov{\vt} 
\intO{ S(\tau_n +, \cdot)} - \omega \ \mbox{for all}\ t \geq \tau_n,\ n \geq n(\omega).
\label{s28}
\end{align}	

\item 

Now, we use the construction of Section \ref{ed}, namely increasing the value of entropy at $\tau_n+$ and concatenating 
the solution as in \eqref{s20}. This procedure yields
\begin{align}
	&\intO{E \Big( (\vr, \vm, \widetilde{S}) (\tau_n +, \cdot ) \Big) } - \Ov{\vt} 
\intO{ \widetilde{S}(\tau_n +, \cdot)}	\br 
&\quad = \mathcal{E}_0 - \Ov{\vt} 
\intO{ \widetilde{S}(\tau_n +, \cdot)} \br
&\quad \leq \intO{E \Big( (\vr, \vm, {S}) (\tau_n +, \cdot ) \Big) } - \Ov{\vt} 
\intO{ {S}(\tau_n +, \cdot)} - G(d(\tau_n)),
\label{s29}	
\end{align}	
cf. \eqref{s20}. 
In addition, 
\begin{align} 
\intO{E \Big( (\vr, \vm, \widetilde{S}) (t\pm, \cdot ) \Big) } &- \Ov{\vt} 
\intO{ \widetilde{S}(t \pm , \cdot)} \leq 
\mathcal{E}_0  - \Ov{\vt} 
\intO{ \widetilde{S}(t \pm , \cdot)} \br
&\leq \mathcal{E}_0 - \Ov{\vt} \intO{ \widetilde{S}(\tau_n +, \cdot)} \ \mbox{for all}\ t > \tau_n.
\label{s30}
\end{align}

Finally, we choose $\omega$ in \eqref{s28} so that 
\begin{equation} \label{s31}
G(d_E(\tau_n)) > 2 \omega \ \mbox{whenever}\ n \geq n(\omega). 
\end{equation}
Combining \eqref{s28}--\eqref{s31} we conclude 
\[
\mathcal{F}(\vr, \vm, \widetilde{S}) < 
\mathcal{F}(\vr, \vm, S) \ \mbox{for a.a.}\ t > \tau_n ,\ n = n(\omega)
\]
in contrast with the assumption that $(\vr, \vm, S)$ is the global minimizer of \eqref{s21}.

\end{enumerate}	
	
We have shown \eqref{s22}; whence completed the proof of Proposition \ref{sP1}.
	
\section{Statistical solutions - the main result}
	\label{m}
	
In the previous part, we have constructed a Borel measurable semigroup 
\begin{equation} \label{m1}
\mathcal{S}: \Big[ t, (\vr_0, \vm_0, S_0, \mathcal{E}_0) \Big] \in [0, \infty) \times \mathcal{P} \mapsto  
\Big( \vr(t, \cdot), \vm(t, \cdot), S(t-, \cdot), \mathcal{E}_0 \Big) \in \mathcal{P}
\end{equation}
consisting of admissible dissipative solutions of the Euler system. 

Suppose we are given a complete Borel probability measure $\sigma \in \mathfrak{P} \Big[ W^{-\ell,2}(\Omega; R^d) \times R \Big]$, 
\begin{equation} \label{m2}
\sigma[ \mathcal{P} ] = 1.
\end{equation}
Next, we associate to $\sigma$ its push forward measure $\sigma_t$ with respect to the semigroup $\mathcal{S}$, 
\begin{equation} \label{m3}
\int_{W^{-\ell,2}(\Omega; R^{d+2}) \times R } \mathcal{G} (\vr_0, \vm_0, S_0, \mathcal{E}_0) \ \D \sigma_t = 
\int_{\mathcal{P}} \mathcal{G} \left( \mathcal{S} \Big[ t, (\vr_0, \vm_0, S_0, \mathcal{E}_0) \Big] \right) \D \sigma \ \mbox{for any}\ t \geq 0,
\end{equation}	
and any $\mathcal{G} \in BC(W^{-\ell,2}(\Omega; R^d) \times R)$.

\begin{Definition}[\bf Statistical solution] \label{mD1}

The family of Markov operators 
\begin{equation} \label{m4}
\mathcal{M}_t : \sigma \in \mathfrak{P}[ \mathcal{P} ] \to \sigma_t \in \mathfrak{P} [ \mathcal{P} ],\ t \geq 0, 
\end{equation}	
defined through \eqref{m3} is called \emph{statistical solution} of the Euler system.	
	
\end{Definition}	

Let us summarize the basic properties of the statistical solutions:

\begin{enumerate}
\item {\bf Markov semigroup.}
The mapping 
\[
\mathcal{M}_t : \mathfrak{P}[\mathcal{P}] \to \mathfrak{P}[\mathcal{P}] 
\]
is a semigroup, meaning 
\begin{equation} \label{m5}
\mathcal{M}_0 = \mathbb{I},\ \mathcal{M}_{t + s} = \mathcal{M}_t \circ \mathcal{M}_s ,\ s,t \geq 0. 
\end{equation}
Moreover, 
\begin{equation} \label{m6}
\mathcal{M}_t \Big[ \delta_{(\vr_0, \vm_0, S_0, \mathcal{E}_0)} \Big] = 
\delta_{\mathcal{S}[t,(\vr_0, \vm_0, S_0, \mathcal{E}_0)]},\ t \geq 0;
\end{equation}
\begin{equation} \label{m7}
\mathcal{M}_t \left[ \sum_{i=1}^n \lambda_i \sigma_i \right] = \sum_{i=1}^n \lambda_i \mathcal{M}_t [\sigma_i]
\ \mbox{whenever}\ \sum_{i=1}^n \lambda_i = 1,\ \lambda_i \geq 0;	
\end{equation}	
\begin{equation} \label{m8}
	t \in [0, \infty) \mapsto \mathcal{M}_t [\sigma ] \ \mbox{is c\` agl\` ad in the narrow topology on}\ \mathfrak{\mathcal{P}}.
\end{equation} 

\item {\bf Regularity.} The semigroup $(\mathcal{M}_t)_{t \geq 0}$ is regular, with the corresponding dual operators

\begin{equation} \label{m9}
\mathcal{M}^*_t [ \mathcal{G}] = \mathcal{G} \circ \mathcal{S}[ t, \cdot]	
\end{equation}	
for any $\mathcal{G} \in BC(\mathcal{P})$, where $\mathcal{G} \circ \mathcal{S}[ t, \cdot]$ is a bounded Borel function on $\mathcal{P}$.	

\item {\bf Vanishing energy defect.}

We define the energy defect as a concave u.s.c function on $\mathcal{P}$, 
\[
d_E = \mathcal{E}_0 - \intO{ E(\vr, \vm, S) } \geq 0.
\]
There holds 
\begin{equation} \label{m10}
\int_{\mathcal{P}} d_E  \ \D \mathcal{M}_t[ \sigma ] \to 0 \ \mbox{as}\ t \to \infty 
\end{equation}	
for any $\sigma \in \mathfrak{P}[\mathcal{P}]$.

\end{enumerate}

\section{Proof of the properties of the statistical solution}
\label{p}	

Our ultimate task is to verify the basic properties of statistical solutions listed in the previous section. 

\begin{enumerate}
	\item The semigroup property \eqref{m5} as well as \eqref{m6}, \eqref{m7} follow immediately from the 
	the semigroup property of the admissible dissipative solutions. 
	
	The c\` agl\` ad continuity stated in \eqref{m8} follows from the analogous property of the 
	solution semigroup, more precisely its entropy component,  
	\[
	\mathcal{S}\left[t, (\vr_0, \vm_0, S_0, \mathcal{E}_0) \right] = (\vr(t, \cdot), \vm (t, \cdot), S(t-, \cdot), \mathcal{E}_0),\ t \geq 0.
	\]
	
	\item The existence of the adjoint semigroup \eqref{m9} follows directly from formula \eqref{m3}.
	
	\item As for the vanishing energy defect claimed in \eqref{m10}, we first rewrite 
	\[
	\int_{\mathcal{P}} d_E \ \D \mathcal{M}_t (\sigma) = \int_{\mathcal{P}} \left[ \mathcal{E}_0 - \intO{ E \Big( (\vr, \vm, S)(t+, \cdot) \Big) } \right] \D \sigma.
	\]
	In accordance with Proposition \ref{sP1}, 
	\[
	\left[ \mathcal{E}_0 - \intO{ E \Big( (\vr, \vm, S)(t+, \cdot) \Big) } \right] \to 0 \ \mbox{as}\ t \to \infty \ \sigma - \mbox{a.s.},
	\]
	and the desired conclusion follows from the Lebesgue Dominance Convergence theorem.
\end{enumerate}	

\section{Conclusion}
\label{c}

We have proposed a construction of statistical solutions to the Euler system of gas dynamics applying a single step selection criterion 
in the class of dissipative solutions. The selection criterion is based on minimizing the Bregman divergence to the maximal entropy 
equilibrium solution in an an exponentially weighted time integral norm. We have shown that this criterion is physically relevant as its 
energy defect decays to zero with growing time. This observation seems consistent with the decay of turbulent phenomena initiated in the flow by the 
data. Indeed the amplitude of the energy defect is directly related to that of Reynolds stress arising in models of deterministic turbulence, 
cf. Bardos et al. \cite{BaGhKa}, \cite{BarNgu}, \cite{BarTi2013}, \cite{BaTi},

An interesting and largely open question is whether or not the selected solution is in fact a weak solution to the Euler system, or, al least, 
its variant defined in \eqref{ss8}. A partial answer is conditional based on a modification of the selection functional 
\[
\int_0^\infty \exp(-\lambda t) \mathcal{F}(\vr, \vm, S)(t, \cdot) \dt,\ \lambda > 0
\]
in the spirit of \cite{FeLmJu2025}. Using the arguments of Section \ref{ed} combined with those of \cite[Section 5]{FeLmJu2025} one could show the following 
result: 

\bigskip

\noindent \textit{Suppose that} 
\[
{\rm arg min} \int_0^\infty \exp(-\lambda t) \mathcal{F}(\vr, \vm, S)(t, \cdot) \dt 
\]
\textit{is independent of $\lambda$ for all $\lambda \geq \lambda_0$. Then the minimizer is a weak solution of the Euler system in the sense of \eqref{ss8}.}

\bigskip
\noindent
We leave the proof of this statement to the interested reader. The resulting solution should be compared to the \emph{absolute energy minimizer} introduced
in the context of the barotropic Euler system in \cite{FeLmJu2025}.

Finally, we compare the results of the present paper with the two step selection procedure proposed in \cite{FeiLM2025II} based on the successive maximization of the entropy functional
\[
\int_0^\infty \exp(-t) \left( \intO{ S (t, \cdot) } \right), 
\]
and the subsequent minimization of the total energy 
\[
\int_0^\infty \exp(-t) \left( \intO{ E\Big( (\vr, \vm, S) (t, \cdot) \Big) } \right).
\]

Since the weak solution are characterized by the vanishing energy defect $d_E$, we easily deduce the following conclusion. 

\bigskip

\noindent 
\textit{Suppose that the admissible dissipative solution introduced in of Definition \ref{sD1} is a weak solution of the 
	Euler system in the sense of \eqref{ss8}. Then it is selected in the first step of the procedure of \cite{FeiLM2025II}. }


\def\cprime{$'$} \def\ocirc#1{\ifmmode\setbox0=\hbox{$#1$}\dimen0=\ht0
	\advance\dimen0 by1pt\rlap{\hbox to\wd0{\hss\raise\dimen0
			\hbox{\hskip.2em$\scriptscriptstyle\circ$}\hss}}#1\else {\accent"17 #1}\fi}

\end{document}